\documentclass{amsart}
\usepackage{graphicx}
\usepackage{color}

\numberwithin{equation}{section}

\theoremstyle{plain}
\newtheorem{theorem}{Theorem}
\newtheorem{lemma}[theorem]{Lemma}
\newtheorem{proposition}[theorem]{Proposition}
\theoremstyle{definition}
\newtheorem{definition}[theorem]{Definition}
\theoremstyle{remark}
\newtheorem{remark}[theorem]{Remark}
\newtheorem*{acknowledgments}{Acknowledgments}

\newcommand\modZ {{\mathbb Z}}
\newcommand\modL {{\mathcal L}}

\begin{document}

\title[Brunnian links, claspers and Goussarov-Vassiliev
  invariants]{Brunnian links, claspers and Goussarov-Vassiliev finite
  type invariants}

\author{Kazuo Habiro}

\address{Research Institute for Mathematical Sciences\\Kyoto
  University\\Kyoto\\606-8502\\Japan}

\email{habiro@kurims.kyoto-u.ac.jp}

\dedicatory{Dedicated to the memory of Mikhail Goussarov}

\begin{abstract}
  We prove that if $n\ge 1$, then an $(n+1)$-component Brunnian link $L$
  in a connected, oriented $3$-manifold is $C_n$-equivalent to an
  unlink.  We also prove that if $n\ge 2$, then $L$ can not be
  distinguished from an unlink by any Goussarov-Vassiliev finite type
  invariant of degree\;$<2n$.
\end{abstract}

\date{October 21, 2005}

\thanks{This research was partially supported by the Japan Society for
the Promotion of Science, Grant-in-Aid for Young Scientists (B),
16740033.}
\maketitle

\section{Introduction}
\label{sec:introduction-1}
Goussarov \cite{Gusarov:91-1,Gusarov:91-2} and Vassiliev
\cite{Vassiliev} independently introduced the notion of finite type
invariants of knots, which provides a beautiful,
unifying view over the quantum link invariants
\cite{Birman:93,Birman-Lin,Kontsevich,Bar-Natan}.
For each oriented, connected $3$-manifold $M$, there is a filtration
\begin{equation*}
  \modZ \modL  = J_0\supset J_1\supset \cdots
\end{equation*}
of the free abelian group $\modZ \modL $ generated by the set $\modL =\modL (M)$ of
ambient isotopy classes of oriented, ordered links in $M$, where for
$n\ge 0$, the subgroup $J_n$ is generated by all the $n$-fold
alternating sums of links defined by `singular links' with $n$ double
points.  An abelian-group-valued link invariant is said to be of
degree$\le n$ if it vanishes on $J_{n+1}$.

\begin{figure}\begin{center}\input{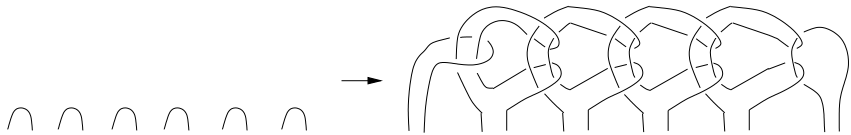}
    \end{center}\caption{A special $n$-variation or a special $C_n$-move ($n=5$).}\label{F02}\end{figure}

\begin{figure}\begin{center}\input{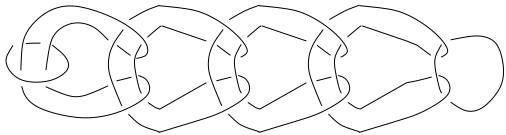}
    \end{center}\caption{Milnor's link of $6$-components.}\label{F01}\end{figure}

Goussarov \cite{Goussarov:finite,Gusarov:variations} and the author
\cite{H} independently introduced theories of surgery along embedded
graphs in $3$-manifolds, which are called {\em $Y$-graphs} or {\em
variation axes} by Goussarov, and {\em claspers} by the author.  For
links, one has the notion of {\em $n$-variation equivalence} (simply
called {\em $n$-equivalence} in \cite{Gusarov:variations}) or {\em
$C_n$-equivalence}, which is generated by {\em $n$-variation}
\cite{Gusarov:variations} or {\em $C_n$-moves} \cite{H}, respectively.
As proved by Goussarov \cite[Theorem 9.3]{Gusarov:variations}, for
string links and knots in $S^3$, the $n$-variation (or $C_n$-)
equivalence is the same as the Goussarov-Ohyama $n$-equivalence
\cite{Gusarov:91-2,Ohyama}.  The $C_n$-equivalence is generated by the
local move depicted in Figure \ref{F02}, i.e., band-summing Milnor's
link of $(n+1)$-components \cite[Figure 7]{Milnor}, see
Figure~\ref{F01}.

One of the main achievements of these theories is the following
characterization of the topological information carried by
Goussarov-Vassiliev finite type invariants.


\begin{theorem}[\cite{Gusarov:variations,H}]
  \label{r13}
  Two knots $K$ and $K'$ in $S^3$ are $n$-variation (or
  $C_n$-)equivalent if and only if we have $K-K'\in J_n$ (i.e., $K$ and
  $K'$ are not distinguished by any Goussarov-Vassiliev invariants of
  degree$<n$.)
\end{theorem}

The variant of Theorem \ref{r13}, with $n$-variation equivalence
replaced by Goussarov-Ohyama $n$-equivalence, is proved previously by
Goussarov \cite{Gusarov:93}.

In \cite[Proposition 7.4]{H}, we observed that for links in $S^3$ there
is a certain difference between the notion of $C_n$-equivalence and
the notion of the Goussarov-Vassiliev finite type invariants of
degree\;$<n$, i.e., Theorem \ref{r13} does not extend to links in
$S^3$.  More specifically, we showed that if $n\ge 2$, then Milnor's
link $L_{n+1}$ of $(n+1)$-components is $C_n$-equivalent but {\em not}
$C_{n+1}$-equivalent to the unlink $U_{n+1}$, but we have
$L_{n+1}-U_{n+1}\in J_{2n}$.  (For $2$-component links, one can easily
observe a similar facts for the Whitehead link $W_2$: $W_2$ is
$C_2$- but not $C_3$-equivalent to the unlink $U_2$, but we have
$W_2-U_2\in J_3$, $\not\in J_4$.)

Note that Milnor's links are examples of {\em Brunnian links}.  Here,
a link $L$ is {\em Brunnian} if any proper sublink of $L$ is an
unlink.  The purpose of this paper is to prove the following results,
which are generalizations of the above-mentioned facts about Milnor's
links to Brunnian links.

Let $M$ be a connected, oriented $3$-manifold.

\begin{theorem}[{Announced in \cite[Remark 7.5]{H} for $M=S^3$}]
  \label{r10}
  For $n\ge 1$, every $(n+1)$-component Brunnian link in $M$ is
  $C_n$-equivalent to an unlink.
\end{theorem}

\begin{theorem}[{Announced in \cite[Remark 7.5]{H} for $M=S^3$}]
  \label{r3}
  Let $n\ge 2$, and let $U$ denote $(n+1)$-component unlink in $M$.  For
  every $(n+1)$-component Brunnian link $L$ in $M$, we have
  $L-U\in J_{2n}$.  (Consequently, $L$ and $U$ can not be distinguished
  by any Goussarov-Vassiliev invariant of degree\;$<2n$ with values in
  any abelian group.)
\end{theorem}

We remark that Theorem \ref{r10} follows from a stronger, but more
technically stated, result (see Theorem \ref{r9} below), which is
proved also by Miyazawa and Yasuhara \cite{Miyazawa-Yasuhara} for
$M=S^3$, independently to the present paper.

\section{Preliminaries}
\label{sec:preliminaries-1}

\subsection{Preliminaries}
\label{sec:preliminaries}
In the rest of this paper, we freely use the definitions, notations
and conventions in \cite{H}.

Throughout the paper, let $M$ denote a connected, oriented
$3$-manifold (possibly noncompact, possibly with boundary).

By a {\em tangle} $\gamma $ in $M$, we mean a proper embedding
$\gamma \colon\thinspace\alpha \rightarrow M$ of a compact, oriented $1$-manifold $\alpha $ into
$M$.  By a {\em link}, we mean a tangle consisting only of circle
components.  (In \cite{H}, tangles are called `links'.)  We sometimes
confuse $\gamma $ and the image $\gamma (\alpha )\subset M$.

Two tangles $\gamma $ and $\gamma '$ in $M$ are {\em equivalent}, denoted by
$\gamma \cong\gamma '$, if $\gamma $ and $\gamma '$ are ambient isotopic fixing the
endpoints.

\subsection{Claspers and tree claspers}
\label{sec:clasp-tree-clasp}
Here we recall some definition of claspers and tree claspers.  See
\cite[\S2, \S3]{H} for the details.

A {\em clasper} $C$ for a tangle $\gamma $ in a $3$-manifold $M$ is a
(possibly unorientable) compact surface $C$ in $\int M$ with some
structure.  $G$ is decomposed into finitely many subsurfaces called
{\em edges}, {\em leaves}, {\em disk-leaves}, {\em nodes} and {\em
boxes}.
We do not repeat here all the rules that should be satisfied by the
subsurfaces.  For the details, see \cite[Definition 2.5]{H}.  We
follow the drawing convention for claspers \cite[Convention 2.6]{H},
in which we draw an edge as a line instead of a band.

Given a clasper $C$, there is defined a way to associate a framed link
$L_C$, see \cite[\S2.2]{H}.  {\em Surgery along $C$} is defined to be
surgery along $L_C$.  A clasper $C$ is called {\em tame} if surgery
along $C$ preserves the homeomorphism type of a regular neighborhood
of $C$ relative to the boundary.  All the clasper which appear in the
present paper are tame, and thus surgery along a clasper can be
regarded as a move of tangle in a fixed $3$-manifold.  The result from
a tangle $\gamma $ of surgery along a clasper $C$ is denoted by $\gamma ^C$.

A {\em strict tree clasper} $T$ is a simply-connected clasper $T$
consisting only of disk-leaves, nodes and edges.  The degree of $T$ is
defined to be the number of nodes plus $1$, which is equal to the
number of disk-leaves minus $1$.  For $n\ge 1$, a {\em $C_n$-tree} will
mean a strict tree clasper of degree $n$.  A {\em $C_n$-move} is
surgery along a $C_n$-tree, which may be regarded as a local move of
tangle since the regular neighborhood of $T$ is a $3$-ball.  The {\em
$C_n$-equivalence} of tangles is the equivalence relation generated by
$C_n$-moves and equivalence of tangles.

A disk-leaf in a clasper is said to be {\em simple} if it intersects
the tangle by one point.  A strict tree clasper is {\em simple} if all
its leaves are simple.

A {\em forest} $F$ will mean `strict forest clasper' in the sense of
\cite[Definition 3.2]{H}, i.e., a clasper consisting of finitely many
disjoint strict tree claspers.  $F$ is said to be simple if all the
components of $F$ are simple.  A {\em $C_n$-forest} is a forest
consisting only of $C_n$-trees.

\section{Brunnian links and $C^a_n$-moves}
\label{sec:brunnian-links-cn}

\subsection{Definition of $C^a_n$-moves}
\label{sec:can-moves}

\begin{definition}
  \label{r7}
  For $k\ge 1$, a {\em $C^a_k$-tree} for a tangle $\gamma $ in $M$ is a
  $C_k$-tree $T$ for $\gamma $ in $M$, such that
  \begin{enumerate}
  \item for each disk-leaf $A$ of $T$, all the strands intersecting
    $A$ are contained in one component of $\gamma $, and
  \item each component of $\gamma $ intersects at least one disk-leaf of
  $T$.  (In other words, $T$ intersects {\em all} the components of
  $\gamma $; this explains `$a$' in `$C^a_k$'.)
  \end{enumerate}
  Note that such a tree exists only when $k\ge l-1$, where $l$ is the
  number of components in $\gamma $.  Note also that the condition (1) is
  vacuous if $T$ is simple.

  A {\em $C^a_k$-move} on a link is surgery along a $C^a_k$-tree.  The
  {\em $C^a_k$-equivalence} is the equivalence relation on tangles
  generated by $C^a_k$-moves.  A {\em $C^a_k$-forest} is a forest
  consisting only of $C^a_k$-trees.
\end{definition}

What makes the notion of $C^a_k$-move useful in the study of Brunnian
links is the following.

\begin{proposition}
  \label{r16}
  A $C^a_k$-move on a tangle preserves the types of the proper
  subtangles.  In particular, if a link $L'$ is $C^a_k$-equivalent to
  a Brunnian link $L$, then $L'$ also is a Brunnian link.
\end{proposition}

\begin{proof}
  Let $T$ be a $C^a_k$-tree for a tangle $\gamma $.  For any proper
  subtangle $\gamma '$, $T$ viewed as a clasper for $\gamma '$ has at least one
  disk-leaf which intersects no components of $\gamma '$.  Hence, by
  \cite[Proposition 3.4]{H}, we have $\gamma '_T\cong \gamma '$.
\end{proof}

Obviously, $C^a_k$-equivalence implies $C_k$-equivalence.  But the
converse does not hold in general, since a $C_k$-move can transform an
unlink into a non-Brunnian link (e.g., a link with a knotted
component).

The following result gives a characterization of Brunnian
links in terms of clasper moves.

Theorem \ref{r10} follows from Theorem \ref{r9} below.

\begin{theorem}
  \label{r9}
  An $(n+1)$-component link $L$ in $M$ ($n\ge 1$) is Brunnian if and
  only if $L$ is $C^a_n$-equivalent to an $n$-component unlink $U$ in
  $M$.
\end{theorem}

As mentioned in the introduction, Theorem \ref{r9} is proved
independently by Miyazawa and Yasuhara \cite{Miyazawa-Yasuhara} for
$M=S^3$.

The rest of this subsection is devoted to proving Theorem \ref{r9}.

The following two lemmas easily follow from the proof of the
corresponding results in \cite{H}.

\begin{lemma}[{$C^a$-version of \cite[Theorem 3.17]{H}}]
  \label{r18}
  For two tangles $\gamma $ and $\gamma '$ in $M$, and an integer $k\ge 1$, the
  following conditions are equivalent.
  \begin{enumerate}
  \item $\gamma $ and $\gamma '$ are $C^a_k$-equivalent.
  \item There is a simple $C^a_k$-forest $F$ for $\gamma $ in $M$ such
    that $\gamma ^F\cong \gamma '$.
  \end{enumerate}
\end{lemma}

\begin{lemma}[{$C^a$-version of \cite[Proposition 4.5]{H}}]
  \label{r4}
  Let $\gamma $ be a tangle in $M$, and let $\gamma _0$ be a component of $\gamma $.
  Let $T_1$ and $T_2$ be $C_k$-trees for a tangle $\gamma $ in $M$,
  differing from each other by a crossing change of an edge with the
  component $\gamma _0$.  Suppose that $T_1$ and $T_2$ are $C^a_k$-trees
  for either $\gamma $ or $\gamma \setminus \gamma _0$.  Then $\gamma ^{T_1}$ and $\gamma ^{T_2}$ are
  related by one $C^a_{k+1}$-move.
\end{lemma}

Now we prove Theorem \ref{r9}.

\begin{proof}[Proof of Theorem \ref{r9}]
  Let $L=L_0\cup L_1\cup \cdots\cup L_n$.

  The `if' part follows since a $C^a_n$-move for an
  $(n+1)$-component link $L$ preserves each proper sublinks of $L$ up
  to isotopy.

  The proof of the `only if' part is by induction on $n$.

  Suppose $n=1$.  Since $L=L_0\cup L_1$ is Brunnian, it follows that both
  $L_0$ and $L_1$ are unknotted in $M$.  In $M$ we can homotop $L_1$
  into an unknot $U_1$, such that $L_0\cup U_1$ is an unlink.  This
  homotopy can be done by ambient isotopy and crossing changes between
  distinct components, i.e., (simple) $C^a_1$-moves.  This shows the assertion.

  Suppose $n>1$.  Since $L$ is Brunnian in $M$, it follows that
  $L'=L\setminus L_0$ is an $n$-component Brunnian link in $M\setminus L_0$.  By
  induction hypothesis, it follows that $L'$ is $C^a_{n-1}$-equivalent
  in $M\setminus L_0$ to an $n$-component unlink $U'$ in $M\setminus L_0$.  By Lemma
  \ref{r18}, there is a $C^a_{n-1}$-forest $F$ for $U'$ in $M\setminus L_0$
  satisfying $(U')^F\cong L'$ in $M\setminus L_0$.  Since $L_0\cup U'$ is an
  unlink, there is a disk $D_0$ in $M$ disjoint from $U'$.
  We may assume that $D_0$ intersects $F$
  only by finitely many transverse intersections with the edges of
  $F$.  By crossing changes between $L_0$ and edges of $F$
  intersecting $D$, we obtain from $L_0$ an unknot $U_0$ in $M$ which
  bounds a disk disjoint from $L'$ and $F$.  By Lemma \ref{r4}, it
  follows that these crossing changes do not change the
  $C^a_n$-equivalence class of the result of surgery.  Hence we have
  \begin{equation*}
    L=L_0\cup L'
    \cong L_0\cup (U')^F
    \cong (L_0\cup U')^F
    \underset{C^a_n}{\sim} (U_0\cup U')^F
    \cong U_0\cup (U')^F
    \cong U_0\cup L'.
  \end{equation*}
  Since $U_0\cup L'$ is an unlink, the assertion follows.
\end{proof}

\subsection{Generalization to tangles}
\label{sec:bottom-tangles}
One can generalize Theorem \ref{r9} to tangles as follows.

Let $c_0,\ldots,c_n\subset \partial M$ be disjoint arcs, and set $c=c_0\cup \cdots\cup c_n$.  A
$(n+1)$-component {\em tangle in $M$ with arc basing $c$} is a tangle
$\gamma $ consisting of $n+1$ properly embedded arcs $\gamma _0,\ldots,\gamma _n$ in $M$
such that $\partial \gamma _i=\partial c_i$ for $i=0,\ldots,n$.  A tangle $\gamma $ with arc
basing $c$ is called {\em trivial} (with respect to $c$) if simple
closed curves $\gamma _i\cup c_i$ for $i=0,\ldots,n$ bounds disjoint disks in $M$.
A tangle $\gamma $ with arc basing $c$ is {\em Brunnian} if every proper
subtangle of $\gamma $ is trivial with respect to the corresponding
$1$-submanifold of $c$.

\begin{theorem}
  \label{r5}
  If $\gamma =\gamma _0\cup \cdots\cup \gamma _n$ ($n\ge 1$) is an $(n+1)$-component tangle in
  $M$ with arc basing $c=c_1\cup \cdots\cup c_n$.  Then $\gamma $ is Brunnian if and
  only if $\gamma $ is $C^a_n$-equivalent to an $(n+1)$-component trivial
  tangle with respect to $c$.
\end{theorem}

\begin{proof}
  Similar to the proof of Theorem \ref{r9}.
\end{proof}

\begin{remark}
  \label{r15}
  The case $M=B^3$ of Theorem \ref{r5} is independently proved by
  Miyazawa and Yasuhara \cite[Proposition 4.1]{Miyazawa-Yasuhara}.
\end{remark}

\begin{remark}
  \label{r26}
  Taniyama \cite{Taniyama} (see also Stanford \cite{Stanford}) proved
  that an $(n+1)$-component Brunnian link is {\em $n$-trivial}, or
  {\em $n$-equivalent} to an unlink.  Here, by `$n$-triviality' and
  `$n$-equivalence' we mean the notion introduced independently by
  Goussarov \cite{Gusarov:91-2} and Ohyama \cite{Ohyama} (see also
  \cite{Taniyama,Gusarov:variations}).  It is well known that
  $C_n$-equivalence implies $n$-equivalence, but the converse seems
  open for links with at least $2$-components.  However, Goussarov
  \cite{Gusarov:variations} proved that $C_n$-equivalence (or
  $n$-variation equivalence) and $n$-equivalence are the same for
  string links in $D^2\times [0,1]$, and hence the case $M=B^3$ of
  Theorem~\ref{r5} follows from the fact (which seems to be well
  known) that $(n+1)$-component Brunnian tangle of arcs in $B^3$ is
  $n$-trivial.
\end{remark}

Using Theorems \ref{r9} and \ref{r5}, we can prove the following fact,
which means that {\em a Brunnian link in $S^3$ is the closure of a
Brunnian tangle in $B^3$.}  (It is clear that, conversely, the closure
of a Brunnian tangle is Brunnian.)

\begin{proposition}
  \label{r2}
  Let $n\ge 2$.  Given an $n$-component Brunnian link $L=L_1\cup \cdots\cup L_n$
  in $S^3$, there is an $n$-component Brunnian tangle
  $\gamma =\gamma _1\cup \cdots\cup \gamma _n$ in a $3$-ball $B^3$ with respect to a basing
  $c=c_1\cup \cdots\cup c_n\subset \partial B^3$ such that the union
  $\bigcup_{i=1}^n\gamma _i\cup c_i\subset B^3\subset S^3$ viewed as a link in $S^3$ is
  equivalent to $L$.
\end{proposition}

\begin{proof}
  By Theorem \ref{r9} and Lemma \ref{r18}, there is a
  simple $C^a_{n-1}$-forest $F$ for an $n$-component unlink
  $U=U_1\cup \cdots\cup U_n$ such that $U^F\cong L$.  Let $D_1,\ldots,D_n$ be
  disjoint discs in $S^3$ bounded by $U_1,\ldots,U_n$, and set
  $D=D_1\cup \cdots\cup D_n$.  Choose a point $p_0\in S^3$ disjoint from $F\cup D$.
  For each $i=1,\ldots,n$, let $p_i\in U_i\setminus F$ and
  let $g_i$ be a simple arc in $M\setminus F$ from $p_0$ to
  $p_i$ such that $g_i\cap D=p_i$.  Here we may assume that
  $g_i\cap g_j=p_0$ if $i\neq j$.
  Let $N$ be a small
  regular neighborhood of $g_1\cup \cdots\cup g_n$, which is a $3$-ball.  Set
  $B^3=\overline{S^3\setminus N}$.  For $i=1,\ldots,n$, set $c_i=\partial B^3\cap D_i$, and
  set $\gamma ^0_i=U_i\cap B^3$.  Then, by Theorem \ref{r5} the result of surgery
  $\gamma =(\gamma ^0_1\cup \cdots\cup \gamma ^0_n)^F$ is Brunnian with respect to
  $c_1\cup \cdots\cup c_n$, and satisfies the assertion.
\end{proof}

\section{Brunnian links and the Goussarov-Vassiliev filtration}
\label{sec:vass-gouss-finite}

\subsection{Definition of the Goussarov-Vassiliev filtration}
\label{sec:defin-gouss-vass}
Here we recall the definition of the Goussarov-Vassiliev filtration for
links using strict tree claspers.  For the details, see \cite[\S 6]{H}.

Let $\modL (M)$ denote the set of equivalence classes of tangles in
$M$.  For $n\ge 0$, define $J_n=J_n(M)\subset \modZ \modL (M)$ as follows.

By a {\em forest scheme} for a tangle $\gamma $ in $M$, we mean a `strict
forest scheme' in the sense of \cite[Definition 6.6]{H}, i.e., a set
$S=\{T_1,\ldots,T_p\}$ of disjoint, strict tree claspers $T_1,\ldots,T_p$ for
a tangle $\gamma $ in $M$.  The {\em degree} of $S$ is defined to be the
sum of the degrees of $T_1,\ldots,T_p$.  Set
\begin{equation*}
  [\gamma ,S]=[\gamma ;T_1,\ldots,T_p]
  =\sum_{S'\subset S}(-1)^{p-|S'|}\gamma ^{\bigcup S'}\in \modZ \modL (M),
\end{equation*}
where the sum is over all subsets $S'$ of $S$, $|S'|$ denotes the
cardinality of $S'$, and $\bigcup S'$ denote the clasper consisting of
the elements of $S'$.

For $n\ge 0$, let $J_n=J_n(M)$ denote the $\modZ $-submodule of $\modZ \modL (M)$
spanned by the elements $[\gamma ,S]$ for any pair $(\gamma ,S)$ of a link $\gamma $
in $M$ and a forest scheme $S$ for $\gamma $ in $M$ of degree $n$.
This defines a descending filtration of $\modZ \modL (M)$:
\begin{equation*}
  \modZ \modL (M)=J_0(M)\supset J_1(M)\supset \cdots,
\end{equation*}
which is the same as the Goussarov-Vassiliev filtration in the usual
sense, defined using singular tangles.

\subsection{Proof of Theorem \ref{r3}}
\label{sec:theorem-proof}

We need some lemmas before proving Theorem \ref{r3}.

\begin{lemma}[{A variant of \cite[Lemma 3.20]{H}}]
  \label{r8}
  Let $\gamma $ be a tangle in $M$, and let $T$ be a strict tree clasper
  for $\gamma $ in $M$.  Let $N$ be a small regular neighborhood of $T$ in
  $M$.  Then the pair $(N,(\gamma \cap N)^T)$ is homeomorphic to
  $(D^2,(\text{$p$ points}))\times [0,1]$, where $p$ is the number of
  points in $T\cap \gamma $.
\end{lemma}

\begin{proof}
  The case where $T$ is simple is a part of \cite[Lemma 3.20]{H}.
  The general case immediately follows from this case.
\end{proof}

\begin{lemma}
  \label{r12}
  Let $1\le n\le r$, and let $L$ be an $(n+1)$-component Brunnian link in
  $M$.  Then there is a forest $F$ for an $(n+1)$-component unlink $U$
  in $M$ satisfying the following properties.
  \begin{enumerate}
  \item $F$ consists of $C^a_l$-trees with $n\le l<r$.
  \item $U$ bounds $n+1$ disjoint disks $D_1\cup \cdots\cup D_{n+1}$ in $M$
  which are disjoint from edges and trivalent vertices of $F$.
  \item $L$ is $C^a_r$-equivalent to $U^F$.
  \end{enumerate}
\end{lemma}

\begin{proof}
  The proof is by induction on $r$.  The case $r=n$ follows
  immediately from Theorem \ref{r9} by setting $F=\emptyset$.  Suppose
  that the result is true for $r\ge n$ and let us verify the case for
  $r+1$.  Let $F$ be as in the statement of the lemma.  Let $N$ be a
  small regular neighborhood of $F$ in $M$.  Then $U^F$ is obtained
  from $U$ by replacing the part $U\cap N$ by $(U\cap N)^F$.  Since $L$ is
  $C^a_r$-equivalent to $U^F$, it follows from Lemma \ref{r18} that
  there is a $C^a_r$-forest $F'$ for $U^F$ such that
  \begin{equation}
    \label{e1}
    (U^F)^{F'}\cong L.
  \end{equation}
  Using Lemma \ref{r8}, we may assume that $F'$ is disjoint from
  $N$, and thus can be regarded as a forest for $U$ disjoint from
  $F$.  Hence we have
  \begin{equation}
    \label{e2}
    (U^F)^{F'}\cong U^{F\cup F'}.
  \end{equation}

  Now $F'$ may intersects $D=D_1\cup \cdots\cup D_{n+1}$.  We may assume that $F'$
  intersects $D$ only by disk-leaves and finitely many transverse
  intersection of $D$ and edges of $F'$.  By Lemma \ref{r4}, without
  changing the result of surgery up to $C^a_{r+1}$-equivalence, we can
  remove the intersection of $D$ and the edges of $F'$ by crossing
  changes between components of $U$ and edges of $F'$ intersecting
  $D$.  Let $F''$ denote the forest obtained from this operation.  Now
  $D$ is disjoint from the edges and trivalent vertices of $F''$, and
  $U^{F\cup F''}$ and $U^{F\cup F'}$ are $C^a_{r+1}$-equivalent.  From this,
  \eqref{e1} and \eqref{e2}, it follows that $F\cup F''$ is a forest with
  the desired properties.
\end{proof}

\begin{definition}
  \label{r14}
  Let $C$ be a clasper for a tangle $\gamma $ in $M$.  We say that a simple
  disk-leaf $A$ of $C$ {\em monopolizes} a circle component $K$ of
  $\gamma $ in $(C,\gamma )$ if there is a $3$-ball $B\subset M$ such that
  $(\gamma \cup C)\cap B$ looks as depicted in Figure \ref{F03}.  We call the
  pair $(A,K)$ a {\em monopoly} in $(C,\gamma )$.  The monopolized
  component $K$ bounds a disk $D$ in $\int M$ which intersect $C$ by
  an arc $A\cap D$.  We call $D$ a {\em monopoly disk} for $K$.
  \begin{figure}\begin{center}\input{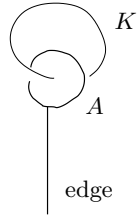}
    \end{center}\caption{A monopoly.}\label{F03}\end{figure}
\end{definition}

\begin{lemma}[Monopoly Lemma]
  \label{r11}
  Suppose $l\ge 1$ and $0\le k\le l+1$ be integers.  Let $T$ be a $C_l$-tree
  for a tangle $\gamma $ in $M$ with $k$ distinct monopolies in $(T,\gamma )$.
  Then we have
  \begin{equation}
  \label{e7}
  \gamma ^T - \gamma  \in  J_{d(l,k)}(M),
  \end{equation}
  where
  \begin{equation*}
    d(l,k) =
    \begin{cases}
      1&\text{if $l=1$, $0\le k\le 2$},\\
      l+k&\text{if $l\ge 2$, $0\le k\le l$},\\
      l+k-1&\text{if $l\ge 2$, $k=l+1$}.
    \end{cases}
  \end{equation*}
\end{lemma}

\begin{proof}
  The case $l=1$ is trivial.  Also, the case $k=l+1$ and $l\ge 2$
  follows from the case $k=l\ge 2$ by ignoring one monopoly.  Hence it
  suffices to prove the case $l\ge 2$, $0\le k\le l$.  Note that if
  $(l,k)=(1,0)$, then we have $d=l+k$.  We will prove by induction on
  $l+k$ that the assertion is true if either $(l,k)=(1,0)$ or $l\ge 2$
  and $0\le k\le l$.

  As we have seen, the case $(l,k)=(1,0)$ is trivial.  Assume
  $l+k\ge 2$.  Let $(A_1,K_1),\ldots,(A_k,K_k)$ be the $k$ monopolies in
  $(T,\gamma )$ with monopoly disks $D_1,\ldots,D_k$, respectively.  Since
  $k\le l$, we can choose one disk-leaf $A_0$ of $T$ distinct from
  $A_1,\ldots,A_k$.  Since $l\ge 2$, $A_0$ is adjacent to a node $Y$.  Let
  $E$ denote the edge between $A_0$ and $Y$.  Let $P'$ and $P''$ be
  the two components of $T\setminus (Y\cup E\cup A_0)$, which are two subtrees in
  $T$.

  Let $l',l''\ge 1$ denote the number of disk-leaves in $P'$ and $P''$,
  respectively.  Let $k'\le l'$ and $k''\le l''$ denote the numbers of the
  monopolizing disk-leaves from $A_1,\ldots,A_k$ contained in $P'$, and
  $P''$, respectively.  We have $l'+l''=l$ and $k'+k''=k$.

  The proof is divided into two cases.

  {\em Case 1.  Either $(l',k')$ or $(l'',k'')$ is $(1,1)$.} We assume
  that $(l',k')=(1,1)$; the other case is proved by the same argument.
  Then $P'$ consists of a monopolizing disk-leaf $A_i$,
  $i\in \{1,\ldots,k\}$, and the incident edge $E'$.  Without loss of
  generality, we may assume that $i=1$.  See Figure \ref{F05} (a).
  \begin{figure}\begin{center}\input{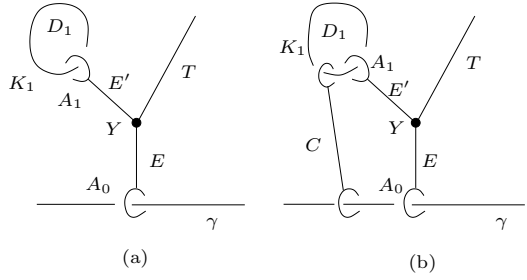}
    \end{center}\caption{Here the lines labeled $\gamma $ depicts a parallel family of
  strands of $\gamma $.}\label{F05}\end{figure} Let $C$ be a $C_1$-tree for $\gamma $ disjoint from
  $T$, as depicted in Figure \ref{F05} (b).  Figure \ref{F06} and
  \cite[Proposition 3.4]{H} imply that $\gamma ^{T\cup C}\cong\gamma ^C$.  (This fact
  is implicit in the proof of \cite[Proposition 7.4]{H}.)
  \begin{figure}\begin{center}\input{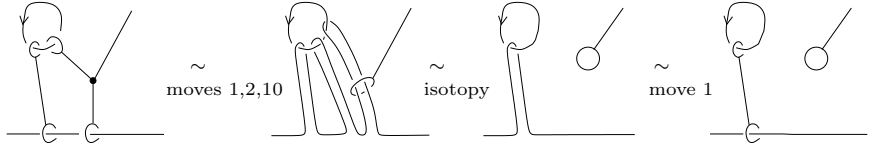}
    \end{center}\caption{Here we use moves 1, 2, 10 of \cite[Proposition 2.7]{H}.
  The orientations given to the circle components are `temporary' and
  may possibly be the opposite to the actual orientations
  simultaneously for all the four figures.}\label{F06}\end{figure} Hence we have
  \begin{equation}
    \label{e8}
    \gamma ^T-\gamma =-(\gamma ^{T\cup C}-\gamma ^T-\gamma ^C+\gamma )=-[\gamma ;T,C].
  \end{equation}
  Let $N$ be a small regular neighborhood of $T\cup D_2\cup D_3\cup \cdots\cup D_k$,
  which is a $3$-ball.  Then, by the induction hypothesis, we have
  \begin{equation}
    \label{e9}
    (\gamma \cap N)^T - \gamma \cap N \in J_{l+k-1}(N).
  \end{equation}
  Since $C$ is a $C_1$-tree, \eqref{e8} and \eqref{e9} implies
  \eqref{e7}.

  {\em Case 2.  Otherwise.}  Apply move 9 in \cite[Proposition 2.9]{H}
  at the disk-leaf $A_0$, see Figure \ref{F04}.  \begin{figure}\begin{center}\input{F04.pstex_t}
    \end{center}\caption{}\label{F04}\end{figure} The result
  is a union $T'\cup T''$ of a $C_{l'}$-tree $T'$ and a $C_{l''}$-tree
  $T''$ for $\gamma $ such that $\gamma ^T\cong \gamma ^{T'\cup T''}$.  Let $N'$ be a
  small regular neighborhood of the union of $T'$ and the monopoly
  disks intersecting $T'$.  Similarly, let $N''$ be a small regular
  neighborhood of the union of $T''$ and the monopoly disks
  intersecting $T''$.  Since $(l',k'),(l'',k'')\neq(1,1)$ and
  $l'+k',l''+k''<l+k$, it follows by induction hypothesis that we have
  \begin{gather*}
    (\gamma \cap N')^{T'}-\gamma \cap N' \in  J_{l'+k'}(N'),\\
    (\gamma \cap N'')^{T''}-\gamma \cap N'' \in  J_{l''+k''}(N'').
  \end{gather*}
  Using \cite[Proposition 3.4]{H}, we see that
  $\gamma ^{T'}\cong\gamma ^{T''}\cong\gamma $.  Hence it follows that
  \begin{equation*}
  \gamma ^T-\gamma =\gamma ^{T'\cup T''}-\gamma ^{T'}-\gamma ^{T''}+\gamma \in J_{l+k}(M)
  \end{equation*}
\end{proof}

\begin{remark}
  \label{r20}
  In \cite[Lemma 7.1]{Conant}, a result similar to Lemma \ref{r11} is
  proved, but it is not strong enough for our purpose.
\end{remark}

Now we prove Theorem \ref{r3}.

\begin{proof}[Proof of Theorem \ref{r3}]
By Theorem \ref{r9} and Lemma \ref{r12} for $r=2n$, there is a forest
$F$ for $U$ in $M$ consisting of simple $C^a_l$-trees with $n\le l<2n$ such
that
\begin{itemize}
\item[(a)] $U$ bounds $n+1$ disjoint disks $D_1,\ldots,D_{n+1}$ in $M$,
disjoint from edges and trivalent vertices of $F$, and
\item[(b)] $L$ is $C^a_{2n}$-equivalent to $U^F$.
\end{itemize}
By the condition (b), we have
\begin{equation}
  \label{e3}
  L- U^F \in  J_{2n}.
\end{equation}
Let $S=\{T_1,\ldots,T_p\}$, $p\ge 0$, be a forest scheme for $U$ in $M$
consisting of the tree claspers $T_1,\ldots,T_p$ contained in $F$.  By an
easy calculation, we have
\begin{equation}
  \label{e5}
  U^F = \sum_{S'\subset S} [U,S'],
\end{equation}
where $S'$ runs over all subsets of $S$.  Since $\deg T_i\ge n$ for all $i$,
we have $\deg S'\ge n|S'|$, where $|S'|$ denotes the number of elements in
$S'$.  Since $|S'|\ge 2$ implies $[U,S']\in J_{2n}$, it follows from
\eqref{e5} that
\begin{equation}
  \label{e4}
  U^F-U \equiv \sum_{i=1}^p [U;T_i] \pmod {J_{2n}}.
\end{equation}
Hence, by \eqref{e3} and \eqref{e4}, it suffices to prove the case
$F=T$ is a $C^a_l$-tree with $n\le l<2n$.  By assumption, there are at
least $k=2n+1-l$ monopolies in $(T,U)$.  Hence by Lemma \ref{r11}, we
have $U^T-U \in J_{d(l,k)}$, where $d(l,k)$ is defined in Lemma
\ref{r11}.  Since $l\ge n\ge 2$, we have $d(l,k)\ge l+k-1\ge 2n$.  Hence we
have $U^T-U\in J_{2n}$.  This completes the proof.
\end{proof}

\subsection{Remarks}
\label{sec:remarks}

\begin{remark}
  \label{r17}
  Przytycki and Taniyama \cite{Przytycki-Taniyama01} proved a
  conjecture by Kanenobu and Miyazawa \cite{Kanenobu-Miyazawa} about
  the homfly polynomial of Brunnian links, and also announced a
  similar result for the Kauffman polynomial.  These results follow
  from Theorem~\ref{r3}.
\end{remark}

\begin{remark}
  \label{r1}
  Yasuhara pointed out to the author that Theorem \ref{r3} implies the
  following generalization.

  {\it Let $n\ge 2$, $m\ge 1$, and let $M$ be a connected, oriented
  $3$-manifold.  Let $L$ and $L'$ be two $(n+1)$-component links in
  $M$ such that
  \begin{enumerate}
  \item both $L$ and $L'$ are $C_m$-equivalent to an $(n+1)$-component
  unlink $U$,
  \item $L$ and $L'$ are $C^a_n$-equivalent to each other.
  \end{enumerate}
  Then we have $L'-L\in J_{l}$, where $l=\min(2n,n+m)$.}

  The proof is as follows.  We may assume that $L=U^F$, where $F$ is a
  $C_m$-forest for $U$.  We may assume also that $L'=U^{F\cup F'}$, where
  $F'$ is a $C^a_n$-forest for $U$, disjoint from $F$.  Then we have
  \begin{equation*}
    L'-L
    =U^{F\cup F'}-U^F
    =(U^{F\cup F'}-U^F-U^{F'}+U)+(U^{F'}-U).
  \end{equation*}
  Here we have $U^{F\cup F'}-U^F-U^{F'}+U\in J_{n+m}$.  We also have
  $U^{F'}-U\in J_{2n}$ by Theorem \ref{r3}.  Hence the assertion.
\end{remark}

\begin{acknowledgments}
  I thank Akira Yasuhara for helpful discussions and comments and for
  asking me about the proof of Theorem \ref{r3} (in the case of
  Brunnian links in $S^3$), which motivated me to write this paper.
  Also, I thank Jean-Baptiste Meilhan for many helpful discussions and
  comments.
\end{acknowledgments}


\begin{thebibliography}{99}

\bibitem{Bar-Natan} D. Bar-Natan, On the Vassiliev knot
  invariants, {\it Topology} {\bf 34} (1995), no. 2, 423--472.

\bibitem{Birman:93} J. S. Birman, New points of view in knot theory,
  {\it Bull. Amer. Math. Soc. (N.S.)} {\bf 28} (1993), no. 2,
  253--287.

\bibitem{Birman-Lin} J. S. Birman and X.-S. Lin, Knot polynomials and
  Vassiliev's invariants, {\it Invent. Math.} {\bf 111} (1993), no. 2,
  225--270.

\bibitem{Conant} J. Conant, {\it Vassiliev invariants and embedded
  gropes}, preprint.

\bibitem{Gusarov:91-1} M. N. Gusarov, A new form of the Conway-Jones
  polynomial of oriented links, (Russian), {\it Zap.
  Nauchn. Sem. Leningrad. Otdel. Mat. Inst. Steklov. (LOMI)} {\bf 193}
  (1991), Geom. i Topol. {\bf 1}, 4--9, 161; translation in {\it
  Topology of manifolds and varieties}, 167--172, Adv. Soviet Math.,
  18, Amer. Math. Soc., Providence, RI, 1994.

\bibitem{Gusarov:91-2} M. Gusarov, On $n$-equivalence of knots and
  invariants of finite degree, {\it Topology of manifolds and
  varieties}, 173--192, Adv. Soviet Math., 18, Amer. Math. Soc.,
  Providence, RI, 1994.

\bibitem{Gusarov:93} M. N. Gusarov, The $n$-equivalence of knots and
  invariants of finite degree, {\it
  Zap. Nauchn. Sem. S.-Petersburg. Otdel. Mat. Inst. Steklov. (POMI)}
  {\bf 208} (1993), 152--173; English transl., {\it J. Math. Sci.}
  {\bf 81} (1996), no. 2, 2549--2561.

\bibitem{Goussarov:finite} M. Goussarov (Gusarov), Finite type invariants and
  $n$-equivalence of $3$-manifolds, {\it C. R. Acad. Sci. Paris
  S\'er. I Math.} {\bf 329} (1999), no. 6, 517--522.

\bibitem{Gusarov:variations} M. N. Gusarov, Variations of knotted
  graphs. The geometric technique of $n$-equivalence.  (Russian), {\it
  Algebra i Analiz} {\bf 12} (2000), no. 4, 79--125; translation in
  {\it St. Petersburg Math. J.} {\bf 12} (2001), no. 4, 569--604.

\bibitem{H} K. Habiro, Claspers and finite type
  invariants of links, {\it Geom. Topol.} {\bf 4} (2000), 1--83.

\bibitem{Kanenobu-Miyazawa} T. Kanenobu and Y. Miyazawa, The second
  and third terms of the HOMFLY polynomial of a link, {\it Kobe
  J. Math.} {\bf 16} (1999), no. 2, 147--159.

\bibitem{Kontsevich} M. Kontsevich, Vassiliev's knot invariants, {\it
  I. M. Gel'fand Seminar}, 137--150, Adv. Soviet Math., 16, Part
  2, Amer. Math. Soc., Providence, RI, 1993.

\bibitem{Milnor} J. Milnor, Link groups, {\it Ann. of Math.} {\bf 59}
  (1954) 177--195.

\bibitem{Miyazawa-Yasuhara} H. A. Miyazawa and A. Yasuhara,
  Classification of $n$-component Brunnian links up to $C_n$-move,
  to appear in {\it Topology Appl.}

\bibitem{Ohyama} Y. Ohyama, A new numerical invariant of knots induced
  from their regular diagrams, {\it Topology Appl.} {\bf 37} (1990),
  no. 3, 249--255.

\bibitem{Przytycki-Taniyama01} J. H. Przytycki and K. Taniyama, The
  Kanenobu-Miyazawa conjecture and the Vassiliev-Gusarov skein modules
  based on mixed crossings, {\it Proc. Amer. Math. Soc.} {\bf 129}
  (2001), no. 9, 2799--2802.

\bibitem{Stanford} T. B. Stanford, Four observations on $n$-triviality
  and Brunnian links, {\it J. Knot Theory Ramifications} {\bf 9}
  (2000), no. 2, 213--219.

\bibitem{Taniyama} K. Taniyama, On similarity of links, {\it Gakujutsu
  Kenkyuu}, School of Education, Waseda University, Series of
  Mathematics 41, 33-36, 1993.

\bibitem{Vassiliev} V. A. Vassiliev, Cohomology of knot spaces, {\it
  Theory of singularities and its applications}, 23--69, Adv. Soviet
  Math., 1, Amer. Math. Soc., Providence, RI, 1990.
\end{thebibliography}
\end{document}